\documentclass{amsart}
\usepackage{amsmath, amsthm, amsopn, amsfonts}
\usepackage{graphicx}
\usepackage{graphpap}

\def\pasdegrille{\let\grille = \pasgrille}

\def\aat#1#2#3{
\divide \dimen1 by 48
\dimen3=\dimen1
\multiply \dimen1 by #1
\advance \dimen1 by -\dimen3
\divide \dimen1 by 101
\multiply \dimen1 by 100
\divide \dimen2 by \count11
\multiply \dimen2 by #2 
\setbox0=\hbox{#3}\ht0=0pt\dp0=0pt
  \rlap{\kern\dimen1 \vbox to0pt{\kern-\dimen2\box0\vss}}\dimen1= \wd1
\dimen2=\ht1}
\def\pasgrille{
\count12= \dimen1 
\divide \count12 by 50
\divide \dimen2 by \count12
\count11 =\dimen2
\ 
\divide \dimen1 by 48
\setlength{\unitlength}{\dimen1}
\smash{\rlap{\ }}
\dimen1= \wd1
\dimen2=\ht1
}
\def\grille{
\count12= \dimen1 
\divide \count12 by 50
\divide \dimen2 by \count12
\count11 =\dimen2
\ 
\divide \dimen1 by 48
\setlength{\unitlength}{\dimen1}
\smash{\rlap{\graphpaper[1](0,0)(50, \count11)}}
\dimen1= \wd1
\dimen2=\ht1
}

\pasdegrille
\usepackage{amssymb}
\usepackage{amsfonts}
\usepackage{graphicx}  
\DeclareGraphicsExtensions{.pdf,.mps}

\usepackage[T1]{fontenc}

\usepackage{epic, eepic}
\def\squarebox#1{\hbox to #1{\hfill\vbox to #1{\vfill}}}

\newcommand{\1}{{\bold 1}}

\newcommand{\CIc}{{\mathcal C}^\infty_{\rm{c}} }

\newcommand{\rest}{\!\!\restriction}

\renewcommand{\Im}{\mathop{\rm Im}\nolimits}

\theoremstyle{plain}

\theoremstyle{definition}

\title
{Bouncing ball modes and quantum chaos}
\author[N. Burq]{Nicolas Burq}
\address{Universit{\'e} Paris Sud,
Math{\'e}matiques,
B{\^a}t 425, 91405 Orsay Cedex}
\email{Nicolas.burq@math.u-psud.fr}
\author[M. Zworski]{Maciej Zworski}
\address{Mathematics Department, University of California \\
Evans Hall, Berkeley, CA 94720, USA}
\email{zworski@math.berkeley.edu}

\usepackage{amssymb}
\usepackage{amsmath, amsthm, amsopn, amsfonts}

\def\Im{\textrm{Im}}
 
\def\11{{\rm 1~\hspace{-1.4ex}l} }

\begin{document}    
\begin{abstract}
Quantum ergodicity of classically chaotic systems has been 
studied extensively both theoretically and experimentally,
in mathematics, and in physics. Despite this long tradition
we are able to present a new rigorous result using only elementary 
calculus. In the case of the famous Bunimovich billiard table shown in 
Fig.1 we prove that the wave functions have to spread into any neighbourhood of the
wings.
\end{abstract}

\maketitle

The quantum/classical correspondence is a puzzling issue that has 
been with us since the advent of quantum mechanics hundred years ago.
Many aspects of it go back to the Newton/Huyghens debate over
the wave vs. corpuscular theories of light.

On the surface of our existence we live in a
world governed by laws of classical physics. That does not
mean that we know precisely how fluids flow or solids move. They are 
described by highly non-linear rules which are hard to unravel mathematically.
Even the simplest classical motion, that of a ball bouncing elastically
from confining walls poses many unanswered questions -- 
see {\small \tt http://www.dynamical\- systems.org/billiard/} for a fun 
introduction. 

If we investigate deeper, or if we simply use any modern technical device,
we come in contact with quantum mechanics. It is governed by a different
set of rules which mix wave and matter. The simplest description
of a wave comes from solving the Helmholtz equation:
\[  ( - \Delta - \lambda^2 ) 
u = 0 \,, \ \ \Delta = \partial_x^2 + \partial_y^2 \,, 
\ \ ( x, y ) \in \Omega \,, \ \ u \rest_{\partial \Omega } = 0 \,. \]
Here we put our wave inside of a two dimensional region $ \Omega $.
In classical wave mechanics the limit $ \lambda
 \rightarrow \infty $ is described
using geometrical optics where the waves propagate along straight lines
reflecting in the boundary $ \partial \Omega $. Roughly speaking, we expect
something similar in the classical/quantum correspondence
with the Helmholtz equation replaced by  its quantum mechanical
version, the Schr\"odinger equation. For many fascinating
illustrations of this we refer to the web art gallery of Rick Heller:
{\small \tt http://www.ericjhellergallery.com}.

\begin{figure}[htb]
\includegraphics[width=9cm]{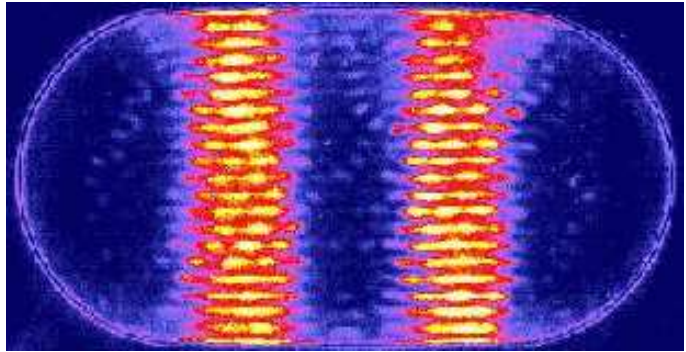}
\caption{An experimental image of bouncing ball modes
in a Bunimovich stadium cavity -- see \cite{ChHu96} and
{\small \tt http://www.bath.ac.uk/$\sim$pyscmd/acoustics}. With 
a certain amount imagination one can see our theorem in this
picture.}
\label{fig:bath}
\end{figure}

Many researchers on different aspects of 
semiclassical analysis have been interested in 
the the correspondence of solutions to the equation above and the 
classical geometry of balls bouncing from the walls of $ \Omega $:
B\"acker, Cvitanovi\'c, 
Eckhardt, Gaspard, Heller, Sridhar, in physics, and Colin de Verdi\`ere,
Melrose, Sj\"ostrand, Zelditch,  in mathematics, to mention some
(see \cite{BaScSt97},\cite{He01}
 for references to the physics literature, and \cite{La93},\cite{MeSj78-1},\cite{Ze03}
for mathematics).
\renewcommand\thefootnote{\dag}

Billiard tables for which the motion is chaotic are a particularly 
interesting model to study\footnote{Of course one would not 
want to play billiards on a table like that, and a completely 
integrable rectangular one can pose enough of a challenge. 
While discussing billiards and the classical/quantum correspondence we 
cannot resist mentioning that 
Pyotr Kapitsa (Nobel Prize in Physics '78) 
was fond of saying that trying
to detect the quantum nature of physical processes at room temperature 
was like trying to investigate the physical laws governing the collision 
of billiard
balls on a table aboard a ship going through rough seas.}. One of the
most famous is the Bunimovich billiard table shown in Fig.1. By adding two 
circular ``wings'' to a rectangular table the motion of a reflecting billiard
ball becomes chaotic, or more precisely, hyperbolic, in the sense that 
changes in initial conditions lead to exponentially large changes in motion as
time goes on. 

As a model for studying quantum phenomena in chaotic systems this billiard
table has become popular in experimental physics. A genuinely quantum example
is shown in Fig.2 -- it comes from the scanning tunnelling microscope work of
Eigler, Crommie, and others~\cite{Cr93}.

\begin{figure}[ht]
\includegraphics[width=10cm]{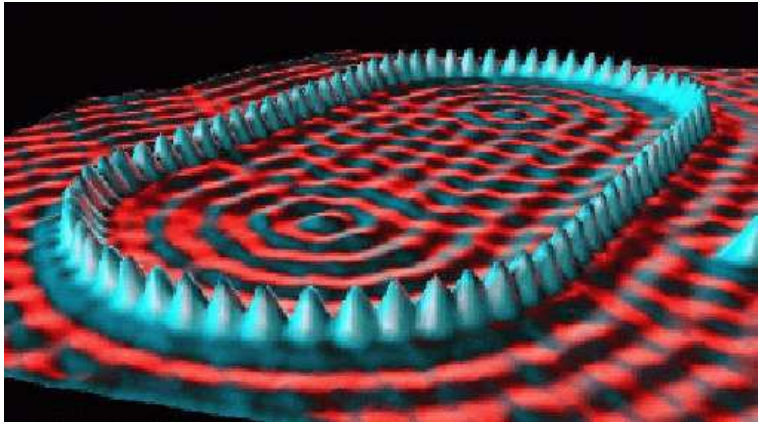}
\caption{Quantum corral in the shape of the Bunimovich stadium. Courtesy of 
IBM{\tiny\textregistered}  Research.}\label{fig:qc}
\end{figure}

One question which is still mysterious to mathematicians and physicists alike
is if the states of this system (that is, solutions of the equation above)
can concentrate on the highly ustable closed orbits of the classical billiard.
{\em Quantum unique ergodicity} states that there is no such concentration
-- see \cite{Sa95},\cite{Sh74},\cite{Ze03} and references given there. In the arithmetic case, that
is for billiards given by arithmetic surfaces where the motion is given by 
the geodesic flow, spectacular advances have been recently achieved by Bourgain, 
Lindenstrauss \cite{BoLi03}, and Sarnak, while for the popular 
quantization of the {\em Arnold cat map} impressive results were
produced by Bonechi, De Bi\`evre, Faure, and Nonnenmacher \cite{FaNo03}, and also by 
Kurlberg and Rudnick \cite{KuRu01}.

Here we describe an elementary but 
striking result in the billiard case.
It follows from adapting the first author's
earlier work in control theory. Although motivated by the more general
aspects of \cite{BuZw03} we give a simple self contained proof.

For a state $ u ( x , y ) $ what counts is its {\em probability density},
$ | u ( x, y ) |^2 dx dy $ -- we assume here that $ |u|^2 
 $ is normalized to have
integral $ 1 $ over $ \Omega $. We say that it is bounded from below in a region if its
integral over that region is bounded from below by a positive constant. 
With this terminology we have, roughly speaking,

\medskip

\noindent
{\bf Theorem.} {\em For any normalized state of the Bunimovich billiard table,
the probability density in any neighbourhood of the
 wings of the table is bounded from below independently
of the energy $ \lambda^2 $.}

\medskip
In particular, the result says that single bouncing ball orbits (that
is orbits following an interval perpendicular to the horizontal straight
boundaries) cannot produce localized waves. Our result allows concentration 
on the full
invariant set of all vertical orbits over $ R $ -- that is 
consistent with the existing physical literature -- 
both numerical and experimental -- see \cite{BaScSt97} and \cite{ChHu96}.
In \cite{BuZw03} we show a stronger result, namely that the 
neighbourhood of the wings can be replaced by any neighbourhood of
the vertical intervals between the wings and the rectangular part.
The proof of that predictable (to experts) improvement is however no longer
elementary and is based on \cite{BaLeRa92}.

The proof of theorem depends on the following unpublished result of the
first author (see \cite{BuZw03} for detailed references and background
material):

\medskip
\noindent
{\bf Proposition.} {\em 
Let $\Delta = \partial_x^2 + \partial_y^2 $, be 
the Laplace operator on the rectangle 
$R= [0, 1]_{x} \times [0,a]_{y}$. Then for any open $\omega\subset R$ 
of the form $ \omega= \omega_{x} \times [0,a]_{y}$ , there 
exists $C$ such that for any solutions of
\[
(-\Delta - \lambda^2) u = f + \partial_x g \ \text{ on $R$}\,, \ \ 
 u \rest_{\partial R}=0 \,, 
\] 
with an arbitrary $ \lambda \geq 0 $
we have 
\[
\int_R | u ( x , y ) |^2 dx dy \leq C \left( 
\int_R ( | f ( x , y )^2 | + | g ( x , y ) |^2 ) dx dy + 
\int_\omega | u ( x , y ) |^2 dx dy \right) \,. 
\]} 
\begin{proof}
We 
decompose $u$, and $ f + \partial_x g $ in terms of 
the basis of $L^2([0,a])$ formed by the Dirichlet eigenfunctions
$e_{k}(y)=  { \sqrt {{2}/a}}\sin(2k\pi y/a)$,
\begin{equation}
u(x,y)= \sum_{k}e_{k}(y) u_{k}(x), \qquad f(x,y) + 
\partial_x g ( x, y ) = \sum_{k}e_{k}(y) ( f_{k}(x) + \partial_x g_k ( x )) 
\end{equation}
we get for $u_{k}, f_{k}$ the equation
\[ 
\left(\partial_{x}^2 + z \right)
u_{k}= f_{k} 
+ \partial_x g_k ,\qquad u_{k}(0)=u_k(1)=0 \,, \ \ 
z = \lambda^2 - \left({2k\pi}/{a}\right)^2 \,.
\] 
It is now easy to see that 
\begin{equation}
\label{eq:easy} 
\int_0^1 |u_{k}( x ) |^2 dx \leq C \left( \int_0^1 ( |f_k ( x ) |^2 + 
| g_k ( x ) |^2 ) dx + \int_{\omega_x } |u_k ( x ) 
|^2 dx \right) \,,   \end{equation}
where $ C $ is independent of $ \lambda_1 $. 
In fact, let us first assume that $ \omega_x = (0, \delta) $, $ \delta > 0 $,
and 
$ z = \lambda_1^2  $, with $ \Im \lambda_1 \leq C $. We then 
choose $ \chi \in \CIc ([ 0,1] ) $ identically zero near $ 0 $ and
identically one on $ [\delta/2,1] $. Then 
\[ \left(\partial_{x}^2 + \lambda_1^2 \right)
(\chi u_{k}) = F_{k} \,, \ \ F_k = \chi ( f_k + \partial_x g_k ) + 
2 \partial_x \chi \partial_x u_k + \partial^2_x \chi u_k \,.\]
We can now use the explicit solution given by 
\[ \chi (x) u_k ( x ) = \frac{1}{\lambda_1} \int_0^x \sin( \lambda_1 ( x - y )) 
F_k ( y ) dy \,. \]
All the terms with $ \partial_x g _k $ and $ \partial_x u _k $ 
can be converted to $ g_k $ and $ u_k $ by integration by parts
(with boundary terms $ 0$ at both ends). Due to the $ \lambda^{-1}_1 $ 
factor that produces no loss and the estimate follows. The argument is 
symmetric under the $ x \mapsto -x $ change, so we can place our control
interval anywhere.

It remains to discuss the case $ z \leq - C < 0 $. Then the estimate
\eqref{eq:easy} follows from integration by parts (where now we do 
not need $ \omega_x $):
\[ \begin{split}
& \int_0^1 \left( f_{k} ( x ) \overline{u_{k} ( x)}  -
 g_{k}  ( x ) \overline{\partial_x u_{k} ( x)} \right) dx = 
 \int_0^1 ( f_{k} ( x ) + \partial_x g_{k} ( x )  ) \overline{u_{k}(x)} dx = \\
& \int_0^1 ( -\partial^2_x - z ) u_{k} ( x) \overline{u_{k} ( x )} dx 
= \int_0^1 \left( |\partial_x u_{k} ( x)|^2 + |z|
| u_{k} ( x)|^2 \right) dx \,. 
\end{split} \]
By the Cauchy-Schwartz inequality, the left hand side is bounded from above by 
\[  \left( \int_0^1 \left( |f_{k} ( x ) |^2 + | g_{k} (  x) |^2 \right)dx 
 \right)^{\frac12} 
\left( \int_0^1 \left( |u_{k} ( x ) |^2  + | \partial_x u_{k}  ( x ) |^2 \right) dx 
 \right)^{\frac12} \,.
\]
Since $ |z| > C > 0 $, \eqref{eq:easy} follows from elementary inequalities
(see \cite[Lemma 4.1]{BuZw03} for a general microlocal argument).
 We can now sum the estimate in $k$ to obtain the proposition.
\end{proof}

We can now present a more precise version of the theorem. For a yet
finer version we refer the reader to \cite[Theorem 3$'$]{BuZw03} and 
\cite[Fig.5]{BuZw03}.

\medskip
\noindent
{\bf Theorem${}'$.}{\em \ Consider $\Omega$ the Bunimovich stadium 
associated to a rectangle $R$. With the convention of Fig.1, let 
$ R_1 $ be any rectangle with the horizontal sides contained in the
sides of $ R $, strictly contained in $R $, and 
with $ R \setminus R_1 $ having two components.

There exists a constant $ C $ depending only on 
 $ \Omega $ and $ R_1 $ such that for any solution of the equation 
\[ (- \Delta - \lambda^2 )v =f\,, \ \  u \rest_{\partial \Omega}= 0 \,,  \ \ 
\lambda \geq 0 \,, \]
we have 
\[ \int_\Omega | v ( x , y ) |^2 dx dy \leq C\left( \int_\Omega|f ( x , y ) |^2 
d x dy + \int_{ \Omega \setminus R_1 }  |v ( x , y ) |^2 d x dy \right) \,.
\]
}

The ``wings'' of the billiard table in the original statement are given by 
$ \Omega \setminus R_1 $. We apply the second theorem with $ f = 0 $ 
to obtain the first one.

\begin{proof}
Let us take $x,y$ as 
the coordinates on the stadium, so that $x$ is the horizontal direction, 
$y$ the vertical direction, and the internal rectangle is 
$[0,1]_{x}\times [0,a]_{y}$. Let us then consider $u$ and $ f$  satisfying
 $(- \Delta- \lambda^2 )u =f$, $ u = 0 $  on the boundary of
the stadium, and $\chi(x)\in \CIc (0,1)$ equal to $1$ on $[\varepsilon, 1- \varepsilon]$. Then $\chi(x) u(x,y)$ is solution of
\[ (- \Delta- \lambda^2 )\chi u =\chi f + [\Delta, \chi] u \ \text{ in $R$}
\] 
with Dirichlet boundary conditions on $\partial R$.
Since $ [\Delta , \chi ] u = 2 \partial_x( \chi' u )  - \chi'' u  $ 
we can apply 
the proposition to obtain 
\[
\int_R |\chi ( x ) u ( x ,y ) |^2  dx dy  \leq 
C \left( \int_R | \chi( x )  f ( x, y ) |^2 dx dy  
+  \int_{ \omega_{\varepsilon} } | u ( x , y) |^2 dx dy \right)
\]
where $\omega_{\varepsilon}$ is a neighbourhood of the support of 
$ \partial_x\chi$. Since we can choose it to be contained in $ R 
\setminus R_1 $, the theorem follows.
\end{proof}

We conclude by remarking that the same argument holds in the
setting discussed recently in \cite{Do03} and \cite{Ze03}, since
in the argument above the rectangle can be replaced by a torus.

\medskip
\noindent
{\sc Acknowledgments.} The authors would like to thank the National 
Science Foundation for partial support under the grant  DMS-0200732.
They are also grateful to Arnd B\"acker, Stephan De Bi\`evre, Richard 
Melrose, 
and Steve Zelditch for helpful 
discussions and comments. They are also grateful to Paul Chinnery,
Victor Humphrey, Don Eigler, and The IBM Corporation for the figures.

\def\cprime{$'$}

\end{document}